\documentclass[a4paper,12pt]{article}
\usepackage[dvips]{color}

\usepackage{amsmath,amssymb}
\usepackage{latexsym}
\usepackage{amsmath}
\usepackage{amssymb}
\usepackage{amsfonts}
\usepackage{epsfig}
\usepackage{amsfonts}
\newtheorem {lemma}{Lemma}

\newtheorem {prop}{Proposition}
\newtheorem {thm}{Theorem}
\def\Eta{\mathfrak{E}}

\def\E{{\mathbb{E\,}}}

\def\P{{\mathbb{P}}}

\newcommand{\Z}{\mathbb{Z}}

\def\|{\,|\,}

\newcommand{\BBox}{\rule{6pt}{6pt}}
\newcommand\Cox{$\hfill \BBox$ \vskip 5mm}

\def \bn{\begin{eqnarray*}}
\def \en{ \end{eqnarray*}}
\def \bnn{\begin{eqnarray}}
\def \enn{ \end{eqnarray}}
\def\M{{\mathcal M}}
\def\m{{{\mathfrak m}}}

\title {Snakes and perturbed random walks}
\author {Gopal Basak\footnote{Statistics and Mathematics Unit, Indian Statistical Institute,
Kolkata 700108, India. Email:~gkb@isical.ac.in}
 \  and Stanislav Volkov\footnote{Department of Mathematics, University of Bristol, BS8~1TW, U.K. Email:~S.Volkov@bristol.ac.uk}}
\begin {document}
\maketitle
\begin{abstract}
In this paper we study some properties of random walks perturbed
at extrema, which are generalizations of the walks considered
e.g.\ in Davis (1999). This process can also be viewed as a
version of {\em excited random walk}, studied recently by many
authors. We obtain a few properties related to the range of the
process with infinite memory. We also  prove the strong law, CLT,
and the criterion for the recurrence of the perturbed walk with
finite memory.
\end{abstract}

\noindent{\bf Key words:} excited / perturbed / cookie random walk, recurrence, transience.

\noindent{\bf AMS 2000 Subject Classification:} Primary 60K35; 60K37.

\section{Introduction}\label{Intro}

This paper has been inspired by the results of Davis (1999) for the random walks perturbed at extrema. Davis (1999) studied the stochastic process which is the limit one of the walks described below. Our purpose here is to study the properties of these perturbed walks as they are, without considering the limit process found in Davis (1999). Let $0<p<1$ and $0<q<1$. Fix a number $L$ which is either a positive integer or $+\infty$. We define a {\em perturbed at $L$-extrema random walk} $X_k$ as a nearest-neighbour random walk on $\Z^1$ with the transitional probabilities $\P(X_{k+1}=x+1\| X_k=x)=1-\P(X_{k+1}=x-1\| X_k=x)$ equal to
 \bn
\left\{
\begin{array}{ll}
p,& \mbox{ if } x=\max_{m=0,1,\dots,L\wedge k} X_{k-m};\\
q,& \mbox{ if } x=\min_{m=0,1,\dots,L\wedge k} X_{k-m};\\
\frac 12 ,& \mbox{ otherwise.}
 \end{array}
\right.
 \en
For definiteness, if $x$ is both the maximum and the minimum, we let this probability be $\frac 12$ (this obviously happens only when $k=0$).

When $L$ is finite, we will call this walk {\em a walk with finite memory}. In this case, it is natural to think of this process as of a ``snake'' (hence the title of the paper, compare our process with a famous video game, released during the mid 1970s, see Surhone, Tennoe, and Henssonow (2010)\,) of length $L$ units moving on the integers, whose transition probabilities depend on whether the snake is surrounded by parts of its body or not.

The case $L=\infty$ corresponds to the walk perturbed at global extrema, and we will refer to this walk as the walk with {\em infinite memory}. It was shown in Davis (1999) that this walk, properly rescaled, converges to a certain stochastic process. Other relevant papers are Davis (1990, 1996); Benjamini and Wilson (2003) and Volkov (2003) studied an excited random walk (ERW), which transitional probabilities differ when the walk visits a site for the first time. Zerner (2005) studied multi-excited random walks on integers. More recently, Basdevant and Singh (2008a) got some interesting results on the speed of ERW, confirming certain conjectures posed in  Zerner (2005), and computed the exact rate of growth in the zero-speed regime in (2008b), while Kosygina and Zerner (2008) obtained annealed CLT for ERW, using branching theory techniques.

When $L=\infty$, our process can be viewed as a special ERW as follows. Consider a site $x$ on a positive axis, and place a geometric number $M_x$ of cookies at the site, such that $\P(M_x=k)=p(1-p)^{k-1}$, $k=1,2,\dots$, and $M_x$, $x=1,2,\dots$, are i.i.d.\ random variables. Whenever the walk visits site $x$ and there are at least two cookie there, it eats one cookie and goes to the left. When it eats the last cookie, it goes to the right. Finally, when there are no cookies left at a site, the walk goes left or right with equal probabilities. On the negative axis, the process is defined similarly with the number of the cookies there distributed according to the law $\P(M_x=k)=(1-q) q^{k-1}$, $x=-1,-2,\dots$. There are no cookies at site $0$.

Note that in Kosygina and Zerner (2008) the number of cookies is uniformly bounded, which is not the case in our model; also our model resembles {\em drilling random walk} introduced in Volkov (2003).

In Section~\ref{infmem}, for $L=\infty$ we focus on the behaviour of the process in {\em finite time}, as opposed to Davis (1999). In particular, we obtain some results on the time it takes before the length of the visited (``cookie-free'') area reaches a certain value. We also compute the limiting probability to be at the global maxima, given that the walk is at one of its extrema.

In Section \ref{finmem}, for $L<\infty$, the process does not resemble ERW mentioned above, and can be viewed as a Markov chain on a product space of $\{-1,+1\}^L \times \Z$. We establish that the finite memory walk is recurrent if and only if $p+q=1$ (compare this with the case $L=\infty$ when the walk is always recurrent). Also we obtain strong law and CLT for our process, and show how the speed of the transient process, when $p+q\ne 1$, decays  as $L$ grows.

We conclude with conjecture and open problems in Section \ref{conj}.

Before we proceed with the next section, let us introduce a few notations. Let
\bn
\M &= \M(k) = \M_L(k)&=\max_{i=1,2,\dots,L\wedge k} X_{k+1-i},
\\
\m &= \m(k) = \m_L(k)&=\min_{i=1,2,\dots,L\wedge k} X_{k+1-i}.
\en
If $X_k\neq \m$ and $X_k\neq\M$ then our process behaves exactly as a simple random walk (SRW). The only differences between $X_k$ and the SRW occur when either $X_k=\m$ or $X_k=\M$. Let ${\cal R}(k)=\M(k)-\m(k)$ be the range (the spread) of the walk and $\rho_n=\min\{k:\ {\cal R}(k)=n\}$ be the stopping time when this range reaches $n$. At some point we will interested in the quantity $\Delta_{n,n+1}=\rho_{n+1}-\rho_n$.

These notations $\M$, $\m$, $\rho_n$, and $\Delta_{n,n+1}$ will be used throughout the paper.

\section{Properties of the perturbed random walks with infinite memory: $L=+\infty$}
\label{infmem}

In this section we establish some interesting properties of the perturbed random walks with infinite memory.

To characterize its distribution, without loss of generality suppose that $\m(k)=0$, $\M(k)=n$, and let $\xi_{l}^{(n)}=\inf\{t:\ X_t=-1\mbox{ or }n+1 \|X_0=l\}$ be the random time until the range is increased given that the walk starts at $l$, $0\le l\le n$. Since $\phi_l=\phi_l^{(n)}(\lambda)=\E \exp(-\lambda\xi_{l}^{(n)})$, the Laplace transform of $\xi_{l}^{(n)}$ satisfies
 \bnn\label{Lapeq}
\phi_0&=&e^{-\lambda}[q\phi_1+(1-q)\cdot 1], \ l=0;\nonumber \\
\phi_l&=&e^{-\lambda}[0.5\phi_{l-1}+0.5\phi_{l+1}], \ 0<l<n;\nonumber \\
\phi_n&=&e^{-\lambda}[p\cdot 1+(1-p)\phi_{n-1}], \ l=n,
 \enn
it can be represented for $\lambda\ge 0$ as
 \bn
\phi_l=C_1\zeta^l+C_2\zeta^{-l},\mbox{ where }
\zeta=e^{\lambda}+\sqrt{e^{2\lambda}-1}.
 \en
The coefficients $C_1,C_2$ should be chosen to satisfy (\ref{Lapeq}). Solving for them, we obtain that
 \bn
\phi_0=
\frac{
(1-q)e^{\lambda}\zeta^2(1-\zeta^{2n})
-(1-p)(1-q)\zeta^3(1-\zeta^{2n-2})+pq\zeta^{n+1}(1-\zeta^{2})
}
{
e^{2\lambda}\zeta^2(1-\zeta^{2n})
-(1+q-p)e^{\lambda}\zeta^3(1-\zeta^{2n-2})
+q(1-p)\zeta^4(1-\zeta^{2n-4})
}
 \en
and
 \bn
\phi_n=
\frac{
p\zeta^2 e^{\lambda}(1-\zeta^{2n})
-pq\zeta^3(1-\zeta^{2n-2})
+(1-p)(1-q)\zeta^{n+1}(1-\zeta^2)
}
{
e^{2\lambda}\zeta^2(1-\zeta^{2n})
-(1+q-p) e^{\lambda}\zeta^3(1-\zeta^{2n-2})
+q(1-p)\zeta^4(1-\zeta^{2n-4})
}
 \en

Let $D_{\M}=D_{\M}(n)=\E(\Delta_{n,n+1}\| X_{\rho_n}=\M(\rho_n))$
and $D_{\m}=D_{\m}(n)=\E(\Delta_{n,n+1}\| X_{\rho_n}=\m(\rho_n))$.
Then
  \bn
  D_{\M}=-\left.\frac{d\phi_n}{d\lambda}\right|_{\lambda=0}
&=&\frac{(1-p-q+pq)n^2+(p+2q-3pq)n+1-p-q+2pq}{1-q-p+2pq+p(1-q)n}\\
D_{\m}-\left.\frac{d\phi_0}{d\lambda}\right|_{\lambda=0}
&=&\frac{pqn^2+(p+2q-3pq)n+1-p-q+2pq}{1-q-p+2pq+p(1-q)n}.
  \en
Note that for large $n$ we have
  \bnn\label{limd}
 D_{\M}(n)=\frac{1-p}p n+O(1), \ D_{\m}(n)=\frac q{1-q}n+O(1).
  \enn

Next, let $p_{\M}=p_{\M}(n)=\P(X_{\rho_{n+1}}=\M(\rho_{n+1}) \| X_{\rho_n}=\M(\rho_n))$ and $p_{\m}=p_{\m}(n)=\P(X_{\rho_{n+1}}=\m(\rho_{n+1}) \| X_{\rho_n}=\m(\rho_n))$ be the probabilities that once the range has increased at the maximum (resp. minimum) the next increase will take place again at  the maximum (resp. minimum). Then these probabilities satisfy
 \bn
p_{\M}=p+(1-p)\left(\frac{n-1}n p_{\M} +\frac 1n (1-p_{\m})\right),
\\
p_{\m}=1-q+q\left(\frac{n-1}n p_{\m} +\frac 1n (1-p_{\M})\right).
 \en
The solution to this system is
 \bnn\label{PMm}
p_{\M}=p_{\M}(n)=\frac{p(1-q)n+pq}{1-p-q+2pq+p(1-q)n},\nonumber\\
p_{\m}=p_{\m}(n)=\frac{p(1-q)n+pq-p-q+1}{1-p-q+2pq+p(1-q)n}.
 \enn
For large $n$ both probabilities are close to one:
 \bn
p_{\M}(n)=1-\frac{1-p}p \times\frac 1n +O(n^{-2}),\
p_{\m}(n)=1-\frac q{1-q}\times\frac 1n +O(n^{-2}).
 \en

Now consider the induced chain $Y_n$ with $Y_n=1$ if $X_{\rho_n}=\M(\rho_n)$, and $Y_n=0$ if $X_{\rho_n}=\m(\rho_n)$. It is straightforward that $Y_k$ is a time-nonhomogeneous Markov chain. We state, however, that it still has a limiting distribution:
\begin{prop}\label{lom}
 \bn
\lim_{n\to\infty}\P(Y_n=1)=\frac{pq}{1-p-q+2pq}=:\pi_{\M}.
 \en
 \end{prop}
{\sf Proof:}
The matrix of transitional probabilities for the chain $\{Y_n\}$ with the states $\{1,0\}$ is
 \bn
A_n=
\left(
\begin{array}{cc}
p_{\M}(n) & 1-p_{\M}(n)\\
1-p_{\m}(n) & p_{\m}(n)
 \end{array}
\right)
=
\left(
\begin{array}{cc}
1-a/n & a/n\\
b/n   & 1-b/n
 \end{array}
\right)
+O(n^{-2})
 \en
where $a=(1-p)/p$, $b=q/(1-q)$, using the formulas for $p_{\M}(n)$ and $p_{\m}(n)$ given by (\ref{PMm}).

Observe that
 \bn
A_n=B C_n B^{-1}+O(n^{-2})
 \en
where
 \bn
B=\left(
\begin{array}{cc}
1 & -a/b \\
1 & 1
 \end{array}
\right)
 \en
and
 \bn
C=\left(
\begin{array}{cc}
1 & 0 \\
0 & 1-(a+b)/n
 \end{array}
\right)
 \en
Hence, for $m\ge n$,
\begin{equation}\label{AAA}
\begin{array}{rcl}
 A_n A_{n+1} A_{n+2}\dots A_m
 &=&B
\left(
\begin{array}{cc}
1 & 0 \\
0 & \prod_{k=n}^m \left(1-\frac{a+b}k\right)
 \end{array}
\right)
B^{-1}
+O(n^{-1}-m^{-1})\\ \\
&=&
\frac 1{a+b}\left(
\begin{array}{cc}
b+a(n/m)^{a+b} & a-a(n/m)^{a+b} \\
b-b(n/m)^{a+b} & a+b(n/m)^{a+b}
 \end{array}
\right)
+O((m-n)n^{-2})
 \end{array}
 \end{equation}
In particular, if $\pi_n=\P(Y_n=1)=b/(a+b)+\delta_n=\pi_{\M}+\delta_n$, then
 \bn
\pi_{n+1}=\P(Y_{n+1}=1)=\pi_n p_{\M}(n)+(1-\pi_n)(1-p_{\m}(n))
=\pi_{\M}+\delta_{n+1}
 \en
where
 \bnn\label{deln}
\delta_{n+1}=\left(1-\frac{a+b}{n}
+O(n^{-2})\right)\delta_n.
 \enn
Since the sequence $\delta_n$'s satisfies (\ref{deln}) and $a+b>0$, it is easy to see that $\lim_{n\to\infty} \delta_n=0$ and therefore $\lim_{n\to\infty} \pi_n=\pi_{\M}$. (In fact, we can even conclude that $\pi_n=\pi_{\M}+O(1)/n^{a+b}$ ).
\Cox

\begin{thm}
(a) For large $n$
 \bn
\E \Delta_{n,n+1}=\frac{q(1-p)n}{1-p-q+2pq} +o(n).
 \en
(b) Asymptotically,
$$
\E \rho_n=\frac{q(1-p)}{1-p-q+2pq} \times \frac {n^2}2+o(n^2).
$$
 \end{thm}
{\sf Proof:}
We start with part (a). Observe that
  \bn
\E \Delta_{n,n+1}/n&=&
\frac {\E(\Delta_{n,n+1}\| X_{\rho_n}=\M(\rho_n))}n \P(Y_n=1)
\\
&+&\frac {\E(\Delta_{n,n+1}\| X_{\rho_n}=\m(\rho_n))}n \P(Y_n=0)
\\
&=&\frac{D_{\M}(n)}n \pi_n+\frac{D_{\m}(n)}n (1-\pi_n)
\longrightarrow \frac {1-p}p \pi_{\M} +\frac q{1-q} (1-\pi_{\M})
  \en
by Proposition~\ref{lom} and formula~(\ref{limd}). Now part (b) of the corollary immediately follows from part (a) and the fact that  $\rho_n=\sum_{i=0}^{n-1}\Delta_{i,i+1}$. \Cox

\section{Finite memory: $ 0<L < \infty$}
\label{finmem}
Note that if we include the history of the process $X_k$ for the past $L$ steps, it becomes a Markov chain. Formally, let $Y_k=(\eta_k^1,\eta_k^2,\dots,\eta_k^L)$ be the sequence of $-1$'a or $+1$'s of length $L$, with $\eta_k^i:=X_{k-i+1}-X_{k-i}$, $i=1,2,\dots,L$. From $Y_k$'s it is possible to extract the information whether the process $X_k$ hit its local maximum or minimum, as described below. Therefore, the pair $(X_k,Y_k)$ is a Markov chain; moreover $Y_k$'s itself form a Markov chain on the space $\Eta^{(L)}=\{-1,+1\}^L$ of the sequences of plus and minus ones of length $L$.

If $Y_k=\eta=(\eta^1,\eta^2,\dots,\eta^L)$, then let  $S(\eta,j)$ denote $\eta^1+\eta^2+\dots+\eta^j$, $j=1,2,\dots,L$. We say that $\eta$ is a local maximum, if $S(\eta,j)\ge 0$ for $j=1,2,\dots,L$; local minimum if $S(\eta,j)\le 0$ for $j=1,2,\dots,L$; and ``neither'' otherwise. Then $X_k$ is a local maximum (minimum resp.) if and only if $Y_k$ is a local maximum (minimum resp.)

Note that from each of the $2^L$ states of $\Eta^{(L)}$ $Y_k=(\eta^1,\eta^2,\dots,\eta^L)$ can go only to two states: $Y_{k+1}=(\eta_*,\eta^1,\eta^2,\dots,\eta^{L-1})$, where $\eta_*=+1$ or $-1$. Observe also that
 \bn
\P(\eta_*=+1)=1-\P(\eta_*=-1)=\left\{\begin{array}{ll}
  p& \mbox{if $Y$ is a local maximum;}\\
  q& \mbox{if $Y$ is a local minimum;}\\
  \frac 12& \mbox{otherwise.}
   \end{array}
\right.
 \en

Since the space $\Eta^{(L)}$ is finite, and $Y_k$ is obviously irreducible, there exists the limiting occupational measure for $Y_k$ denoted as $\pi=\{\pi^{(L)} (\eta)\}_{\eta\in\Eta^{(L)}}$, which obviously depends on $p$, $q$, and $L$.

Let
  \bn
 \pi(\max)=\sum_{\eta\mbox{ is local max.}} \pi(\eta)
  \en
and
  \bn
 \pi(\min)=\sum_{\eta\mbox{ is local min.}} \pi(\eta)
  \en
Observe that  $0<\pi(\min)<1$ and $0<\pi(\max)<1$, and set $\Delta=\Delta_L:=(2p-1)\pi(\max) + (2q-1)\pi(\min)$.

\begin{lemma}\label{lm:delta}
 \bn
  \lim_{n\to\infty} \frac{X_n}n=\Delta_L \mbox{ a.s.}
 \en
 and hence if $\Delta_L\ne 0$ then $X_n$ is transient.
 \end{lemma}
{\sf Proof of the Lemma} Let $f:\Eta^{(L)}\to \{-1,+1\}$ be  such that $f(\eta)=\eta^1$, i.e.\ the first coordinate of $\eta$. By the strong law for the Markov chains (see e.g.~\cite{BHA}, p.145)
 \bn
 \lim_{n\to\infty} \frac 1n \sum_{m=1}^n f(Y_m)=\E_{\pi} f(Y_1)
 \en
Since by the construction  $Y_k$'s we have $\sum_{m=1}^n f(Y_m)=\sum_{m=1}^n (X_m-X_{m-1})=X_n-X_0$, it suffices to show $\E_{\pi} f(Y_1)=\Delta$.

Indeed,
 \bn
\E_{\pi} f(Y_1)=\E_{\pi} Y_1^1
&=&\pi(\max)\E_{\pi} (Y_1^1\| Y_0 \mbox{ is loc.max.})\\
&+&\pi(\min)\E_{\pi} (Y_1^1\| Y_0 \mbox{ is loc.min.})\\
&+&(1-\pi(\max)-\pi(\max))\E_{\pi} (Y_1^1\| Y_0 \mbox{ is neither})\\
&=&\pi(\max)(2p-1)+\pi(\min)(2q-1)\\
&+&(1-\pi(\max)-\pi(\max))\times 0 =\Delta.
 \en
The last statement of the Lemma is straightforward.
\Cox

\begin{lemma}\label{lm:pq1}
 $\Delta<0$, $=0$, or $>0$   if  $p+q -1<0$, $=0$, or $>0$ respectively.
\end{lemma}
{\sf Proof:}
\underline{Case 1: $p+q=1$.}
In this case, the chain $Y_n$ is symmetric. Indeed, if we replace each $+1$ by $-1$ and vice versa, it will have the same distribution since $q=1-p$. Thus by symmetry $\pi(\max)=\pi(\min)$. whence $\Delta=(2p-1)\pi(\max) + (2q-1)\pi(\min) = ((2p-1) + (2q-1))\pi(\max) = 2(p+q-1)\pi(\max) = 0$.

\underline{Case 2a: $p+q>1$, $p,q \ge \frac 12$ .}
In this case, $2p-1\ge 0$ and $2q-1>0$ or $2p-1> 0$ and $2q-1\ge 0$. Hence, $ \Delta=(2p-1)\pi(\max) + (2q-1)\pi(\min) > 0$

\bigskip
\underline{Case 2b: $p+q>1$, $p \ge \frac 12\ge q$.} Construct a new {\it symmetric} chain $\{\tilde Y_n\}$  with $\tilde p=p$ and $\tilde q=1-\tilde p<q\le \frac 12$. Denote $\pi$ and $\tilde \pi$ the stationary distributions for $Y_n$ and $\tilde Y_n$ respectively. Using a sequence of i.i.d.\  $Uniform[0,1]$ random variables $U_1,U_2,\dots$ we will use couple  $Y_n$ and $\tilde Y_n$ to demonstrate that $\pi(\max)\ge \tilde \pi(\max)$ and $\pi(\min)<\tilde \pi(\min)$. Then the comparison with case 1 will
imply $\Delta>0$.

We say that $Y_n$ {\it goes right} ({\it left} resp.), if in the notations of this section $\eta_*=+1$ ($-1$ resp.). The analogous terminology is used for $\tilde Y_n$.

The rules of the coupling are quite standard and are as follows. If  given $Y_{n-1}=\eta$ the probability to go right is $x$, where $x=p$, $q$, or  $\frac 12$ resp. ($\eta$ is local maximum, or minimum, or neither resp.), then $Y_n$ goes right if $U_n\le x$ and goes left otherwise. Similarly, $\tilde Y_n$ goes to the right from the state $\tilde \eta$ if and only if $U_n\le \tilde x$ where $\tilde x=\frac 12$,  $\tilde q$ or $\tilde p$ depending on the state $\tilde \eta$.

For $\eta,\tilde \eta\in\Eta^{(L)}$ we write $\eta\succeq\tilde \eta$ whenever $\eta^i\ge (\tilde \eta)^i$ for all $i=1,2,\dots,L$. Start with $Y_0=\tilde Y_0$. Then it is easy to see by induction that $Y_n\succeq\tilde  Y_n$ for all $n$. Indeed, $Y_0\succeq \tilde Y_0$. Next, if $Y_m\succeq \tilde Y_m$ then

(a) if $Y_m$ is a local minimum then  $\tilde Y_m$ is also a local minimum and hence $\tilde \eta\le \eta$ with the strict inequality whenever $\tilde q<U_m\le q$;

(b) if $Y_m$ is a local maximum then
$\tilde \eta\le \eta$ since $p\ge \max(\tilde p,\frac 12,\tilde q)$

(c) if $Y_m$ is neither of the above then $\tilde Y_n$ cannot be local maximum, and again $\tilde \eta\le \eta$ since $\frac 12 \ge \max(\frac 12, \tilde q)$.

In all three cases above we have $Y_{m+1}\succeq \tilde Y_{m+1}$.

Finally, ``$Y_n$ is a local minimum'' implies ``$\tilde Y_n$ is a local minimum'' and also ``$\tilde Y_n$ is a local maximum'' implies ``$Y_n$ is in a local maximum'', whence $\pi(\max)\ge \tilde \pi(\max)$ and $\pi(\min)\le \tilde \pi(\min)$. Moreover, the event $\{\tilde q<U_n\le q\}$ has a positive probability, consequently a positive fraction of times $\tilde Y_n$ will achieve new local minimum right after it is a local minimum, while $Y_n$ will go the right. Hence, $\tilde \pi(min)>\pi(\min)$.

\bigskip
\underline{Case 2c: $p+q>1$, $p \le \frac 12 \le q$.} Here we will have to construct a series of couplings, as the argument of Case 2b unfortunately cannot be applied directly.

We will also work directly with $X_n$ rather than with $Y_n$; clearly, $\pi(\max)$ and $\pi(\min)$ are the values determined also by the process $X_n$.

First, we construct the second  process $\{\tilde X_n\}$ similar to $X_n$ with $\tilde q=q$ but with $\tilde p=1-q<p$. Observe that the process $\tilde X_n$ is symmetric and define the corresponding Markov chain $\tilde Y_n$ for $\tilde X_n$ in the same way $Y_n$ was defined. Start with $X_0=0=\tilde X_0$. Draw a $Uniform(0,1)$ random variable and go to the right if this variable  s less than $1/2$ and left otherwise. At the $k$-th stage, draw an independent $Uniform(0,1)$ random variable, and move to the right for the asymmetric chain, $\{X_n,Y_n\}$
$$
\left\{\begin{array}{llll}
\mbox{ if} \ U_k < p & \mbox{ whenever } X_{k-1} \mbox{ is local max;} \\
\mbox{ if} \ U_k < q & \mbox{ whenever } X_{k-1} \mbox{ is local min;}\\
\mbox{ if} \ U_k < \frac 12 & \mbox{ whenever $X_{k-1}$ is neither.}
 \end{array}
\right.
$$
Similarly, move to the right for the symmetric chain, $\{\tilde X_n,\tilde Y_n\}$,
$$
\left\{\begin{array}{llll}
\mbox{ if} \ U_k < \tilde p & \mbox{ whenever }
  \tilde X_{k-1} \mbox{ is local max;} \\
\mbox{ if} \ U_k < \tilde q=q & \mbox{ whenever }
  \tilde X_{k-1} \mbox{ is local min;}\\
\mbox{ if} \ U_k < \frac 12 & \mbox{ whenever $\tilde X_{k-1}$ is neither.}
 \end{array}
\right.
$$
Observe that both processes follow the same path till they hit a local maximum, say, at the $(k-1)$-st step.  Then the symmetric process has the smaller probability $\tilde p$ to move to the right (i.e., another local maxima) than the original asymmetric process. In fact, if they make different moves, then the asymmetric process satisfies $X_k=\tilde X_k+ 2$.

{\sf Our second step} is to show by induction that $X_k\ge \tilde X_k$, in fact, for all $k$. Assume that $\tilde X_k\le X_k$, for all  $k\le m$, and if for some $k$, $\tilde X_k < X_k$, then $\tilde X_k+2 \le X_k$. We now show that asymmetric process goes ahead of the symmetric process.

In the table below we write all the possibilities for the two processes at the $m$-th step, and then in each cell we write the pair of probabilities to move to the right. The first number in the brackets is the probability for $X_n$ and the second number is the one for the symmetric process $\tilde X_n$.
\bigskip

\begin{center}
\begin{tabular}{|l|c|c|c|}
\hline
$X_m\backslash \tilde X_m$  & local max. & neither    & local min. \\
\hline
local max.   & $(p,\tilde p)$   & $(p, \frac 12)$   & $(p, q)$    \\
 neither    & $(\frac 12,\tilde p)$ & $(\frac 12, \frac 12)$ & $(\frac 12, q)$  \\
local min.   & $(q, \tilde p)$   & $(q, \frac 12)$   & $(q, q)$    \\
\hline
 \end{tabular}
\end{center}

\bigskip
In all the cells of lower triangular positions including diagonals (i.e., cells with coordinates
 $(1,1)$, $(2,1)$, $(2,2)$, $(3,1)$,  $(3,2)$,  $(3,3)$\ ) the probability of moving to the right for the $X_n$ is bigger than or equal to the probability of moving to the right for the $\tilde X_n$, hence $X_{m+1} \ge \tilde X_{m+1}$ for these cells.

Now for cells $(1,2)$ or $(1,3)$,  if $X_m$ is at the local maximum and $\tilde X_m=X_m$, then $\tilde X_m$ must also be at the local maximum  (since $\tilde X_k\le X_k, \ \forall k\le m$ by the assumption of induction), which contradicts the fact that $\tilde X_m$ is not at local maximum.  Therefore, $\tilde X_m+2\le X_m$. Hence at the $(m+1)$-st step if the processes do move not in the same direction, then they would be at most equal, i.e., $\tilde X_{m+1} \le X_{m+1}$.

Finally, for cell $(2,3)$, if $\tilde X_m$ is at local minimum and $X_m=\tilde X_m$, then $X_m$ must also be at the local minimum (since $X_k\ge \tilde X_k, \ \forall k\le m$), which contradicts the fact that $X_m$ is at neither local max nor local min. Thus again $\tilde X_m+2\le X_m$ and if they make opposite moves at the $(m+1)$-st step they would be at most equal:  $\tilde X_{m+1} \le X_{m+1}$ as before.

Hence we conclude the proof by induction that $X_k\ge \tilde X_k$ for all $k$ and some times $\tilde X_k+2 \le X_k$.

{\sf The third step} consists in using regeneration arguments. Let $\tau_1$ be the first time $k$ when $Y_k=\tilde Y_k=(+1,+1,\dots,+1)$, which is obviously finite as it is stochastically bounded by a geometric random variable with parameter $\max(p,q,\tilde p,\tilde q, 0.5)$. At time $\tau_1$ construct a copy  $\tilde X_k^{(1)}$ of the process $\tilde X_k^{(0)}:=\tilde X_k$, such that $\tilde Y_{\tau_1}^{(1)}=Y_{\tau_1}$. By the arguments of the second step, $\tilde X_k\ge \tilde X_k^{(1)}$ for all $k\ge \tau_1$. Now let  \bn \tau_2=\min\{k>\tau_1 :\ Y_k=\tilde Y^{(1)}_k=(+1,+1,\dots,+1)\}  \en which is finite by the argument above, and construct another copy   $\tilde X_k^{(2)}$ of the process $\tilde X_k$, now such that $\tilde Y_{\tau_2}^{(1)}=Y_{\tau_2}$. Repeating this procedure indefinitely we construct the sequence of stopping times $\tau_1,\tau_2,\tau_3,\dots$ and the sequence of processes $\tilde X_k^{(1)},\tilde X_k^{(2)},\tilde X_k^{(3)},\dots$. Moreover, $\tau_{n+1}-\tau_n$, $n\ge 1$ are i.i.d.\  with finite expectation, say $\bar\tau$. Also, $\tilde X_k\ge \tilde X_k^{(n)}$ for all $k\ge \tau_n$, and because of the arguments of step two, $\tilde X_{\tau_n}^{(n-1)} - \tilde X_{\tau_n}^{(n)}\ge 2$ with a positive probability; moreover, these increments are nonnegative and  independent for different $n$.

Observing that
 \bn
\sum_{m=1}^n
\left[\tilde X_{\tau_m}^{(m)} - \tilde X_{\tau_m}^{(m-1)}\right]
 =X_{\tau_n}-\tilde X_{\tau_n}
 \en and using the strong law we obtain  \bn
 \frac{\tau_n}n &\to & \bar\tau\hskip 10mm \mbox{ a.s.},\\
 \frac{X_{\tau_n}- \tilde X_{\tau_n}}n & \to & \alpha>0\hskip 5mm \mbox{ a.s.}
 \en
for some constant $\alpha$. This, in turn, implies that
 \bn
 \frac{X_k- \tilde X_k}k\to \alpha/\bar\tau>0
 \en
after applying the renewal law of large numbers (Theorem 1.7.3 from \cite{DUR}) for the sequence $\{\tau_k\}$.

Finally, in our {\sf fourth step}, we use Lemma~\ref{lm:delta} and our result for the symmetric Case 1 applicable to $\tilde X_k$, to conclude that $\Delta>\tilde \Delta=0$. Note that this step is only required for the proof of Lemma~\ref{lm:pq1} and is not needed for the following Theorem~\ref{thfinmem}, as the transience of $X_n$ follows immediately from step 3 above (and, of course, Case 1).

\bigskip
\underline{Case 3a: $p+q <1$, $p,q \le \frac 12$ .} In this case,  $2p-1\le 0$ and $2q-1\le 0$ and at least one of the inequalities is strict. Consequently, $\Delta= (2p-1)\pi(\max) + (2q-1)\pi(\min) < 0$.

\bigskip
\underline{Case 3b: $p+q <1$, $p \le \frac 12\le q$.} The proof is exactly similar to that of Case 2c  with role of $p$ replaced by $1-q$ and $q$ by $1-p$, yielding $\Delta= (2p-1)\pi(\max) + (2q-1)\pi(\min) < 0$.

\bigskip
\underline{Case 3c: $p+q <1$, $p \ge \frac 12\ge q$.} In this case, the proof is exactly similar to that of Case 2b  with role of $p$ replaced by $1-q$ and $q$ by $1-p$, whence $ \Delta=(2p-1)\pi(\max) + (2q-1)\pi(\min) < 0$.
\Cox

We now give the criteria for recurrence / transience of the finite memory chain.

\begin{figure}[htb]
\centerline{\hbox{\psfig{figure=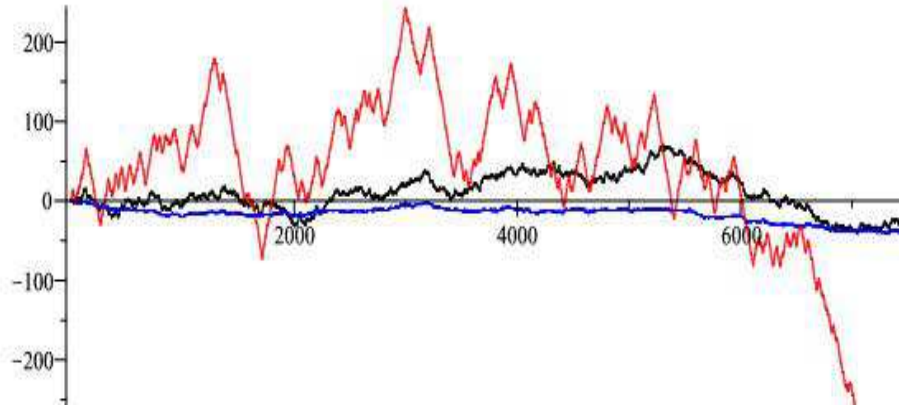,width=12cm,height=6cm}}}
\caption{Trajectories of recurrent walks: red: $p=1-q=0.9$, black: $p=q=1/2$,  blue: $p=1-q=0.1$.}
\label{figSRWvsSNAKE}
\end{figure}

\begin{thm}\label{thfinmem}
For any $L\ge 1$, the finite memory chain is recurrent if and only if $p+q = 1$.
\end{thm}
We must note that even when the chain is recurrent but $p\ne 1/2$, the behaviour of the walk is different from that of a simple random walk. In particular, the variability increases for higher values of $p(=1-q)$. See Figure~\ref{figSRWvsSNAKE}.

\noindent
{\sf Proof of Theorem~\ref{thfinmem}:}
If $p+q\ne 1$ then by Lemma~\ref{lm:pq1} $\Delta\ne 0$ and
hence by Lemma~\ref{lm:delta} $|X_n|\to\infty$ a.s.

If $p+q=1$ then the chain $(X_n,Y_n)$ is symmetric with respect to the change $-1 \leftrightarrow +1$. Let $T_{+\infty}=\{X_n\to+\infty\}$ and $T_{-\infty}=\{X_n\to -\infty\}$. By symmetry, $\P(T_{+\infty})=\P(T_{-\infty})$. On the other hand, $\{X_n\to+\infty\}$ is a tail event, since, for example,  there are infinitely many regeneration times when $X_i=X_{i-1}+1$ for $i=n,n-1,n-2,\dots,n-L+1$. Therefore, by Kolmogorov's zero-one law $\P(T_{+\infty})\in\{0,1\}$. Hence $\P(T_{+\infty})=\P(T_{-\infty})=0$ and $\P(|X_n|\not\to \infty)=1$. Consequently, there is an $k \in\Z$ such that $X_n=k$ for infinitely many $n$. And every time the walk hits $k$, the probability it will reach $0$ in $k$ steps is at least $[\min(p,1-p,q,1-q,\frac 12)]^{|k|}>0$. This implies the recurrence of $X_n$.
\Cox

\begin{thm}\label{thclt}
For any $L\ge 1$, the finite memory chain satisfies the central limit theorem, that is
$$
\frac{X_n-n \Delta_L }{\sqrt{n}}\Longrightarrow {\cal N}(0,\sigma^2)
$$
where $\sigma^2=Var_{\pi}(f(Y_1)) + 2 \lim_{\{m\to \infty\}} \sum_{k=1}^{m+1} Cov_{\pi}(f(Y_1), f(Y_{k+1}))$.
\end{thm}
{\sf Proof:} For $f$ defined as in Lemma \ref{lm:delta}, use Functional CLT on positive recurrent Markov chains $\{Y\}$ (see, Theorem 10.2, p.150 of \cite{BHA}) to get
$$
\frac{X_n-n \Delta_L }{\sqrt{n}} =
\frac {1}{\sqrt{n}} \big [\sum_{m=1}^n (f(Y_m)-\E_{\pi} f(Y_1)) + X_0 \big ]
\Longrightarrow
{\cal N}(0,\sigma^2)
$$
where
$\sigma^2 = Var_{\pi}(f(Y_1)) + 2 \lim_{\{m\to \infty\}} \sum_{k=1}^{m+1} Cov_{\pi}(f(Y_1), f(Y_{k+1}))$
whenever the limit exists and is finite, which holds for a finite state-space Markov Chain $\{Y\}$.
\Cox

\section{Conjectures and open problems}
\label{conj}

Here we list a few open problems and conjectures at which we have arrived by mostly looking at simulations of the process.

\begin{figure}[htb]
\centerline{\hbox{\psfig{figure=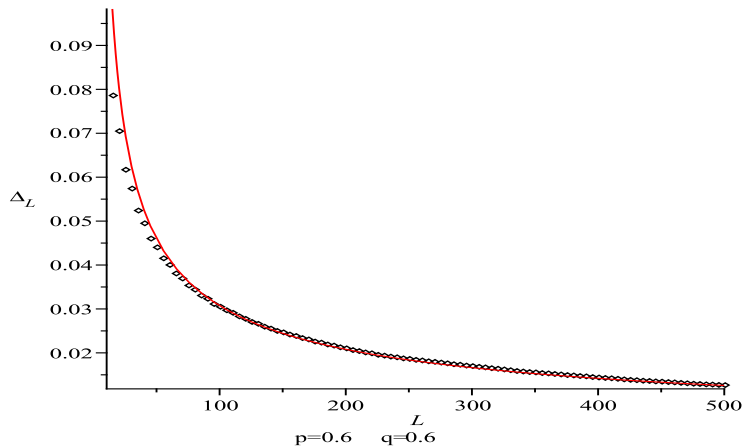,width=10cm,height=6cm}}}
\caption{The speed $\Delta_L$ as a function of $L$ vs. $const\cdot (2L\log\log L)^{-1/2}$.}
\label{FigSpeed}
\end{figure}

In the transient case, when $p+q\ne 1$, the numerical simulations suggest that for fixed $p$ and $q$ we have
$$
\Delta_L \propto \frac 1{\sqrt{2 L\log\log L}}
$$
(see Figure~\ref{FigSpeed}). We believe that this order of magnitude corresponds to the fact that the range of the walk within the last $L$ steps is of order $\sqrt L$ and hence the frequency at which it visits the local maxima and minima, where it gets ``a push'' is something like $L^{-1/2}$ but unfortunately we do not have proof of this fact. The intuition behind this is that for a simple random walk ($p=q=1/2$) the probability to be at the maxima is asymptotically $1/\sqrt{\pi L}$; this follows from Theorem 1.a in Chapter XII.8 and Theorem 1 in Chapter XVIII.5 in~\cite{WF}.

\begin{figure}[htb]\protect
\centerline{\hbox{\psfig{figure=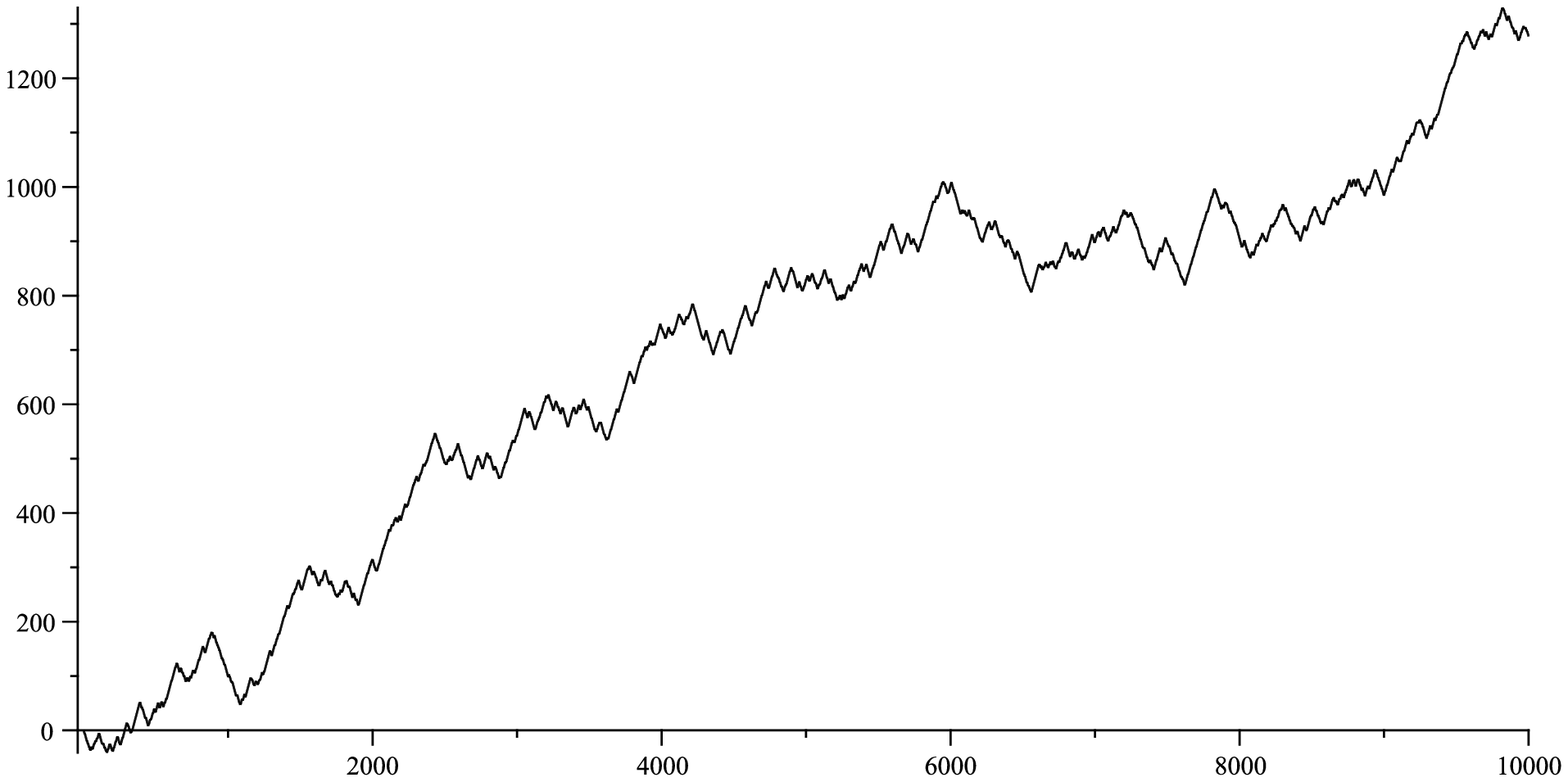,width=8cm,height=6cm}}
            \hbox{\psfig{figure=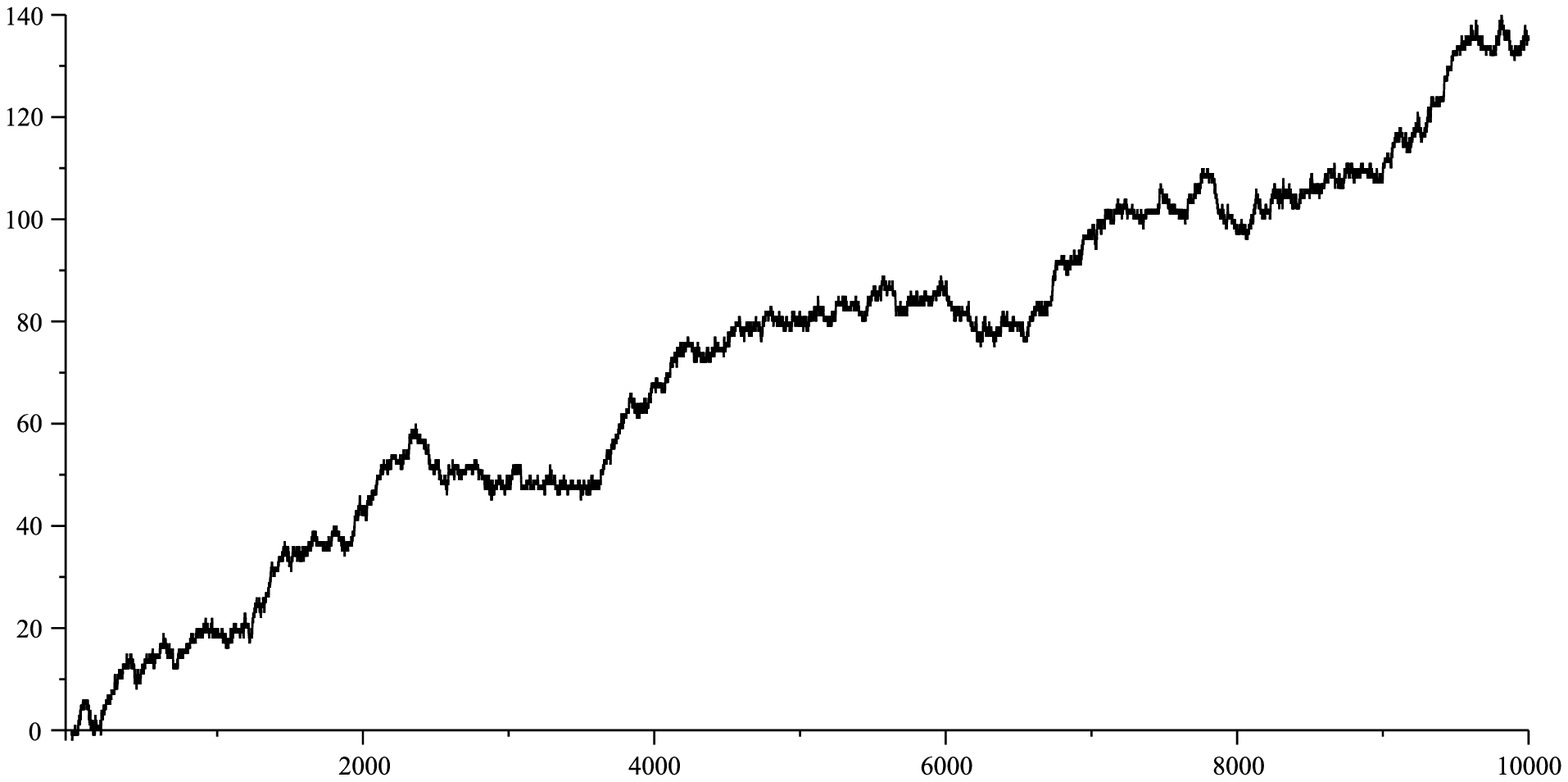,width=8cm,height=6cm}}}
\caption{Left: $p=0.95$, $q=0.15$. Right: $q=0.95$, $p=0.15$.}
\label{FigPq}
\end{figure}

Also, we conjecture that $\Delta_L$ depends not only on ``drift'' $p+q-1$, but in a complicated way on both $p$ and $q$, see Figure~\ref{FigPq} where in both cases the walk is transient.

Recall that in general we have $\Delta_L=(2p-1)\pi(\max) + (2q-1)\pi(\min)$, so estimating $\pi(\max)$ and $\pi(\min)$ is crucial in order to get the speed of the walk. We have another conjecture justified numerically: if $q=1/2$ (so the walk is not perturbed at the minima) then
$$
\pi_{\max}=\frac 1{1+a_L(1-p)},\ \, \ a_L\sim L^{1/2}.
$$
Again, we do not have a rigorous argument for this, and it would be hence nice to obtain a rigorous proof of this asymptotic dependence.


\begin{thebibliography}{99}

\bibitem{Bsingh2}
Basdevant, A-L.,  Singh, A. (2008b) Rate of growth of a transient cookie random walk. {\it Electron. J. Probab.} {\bf 13}, paper~no.~26, 811--851.

\bibitem{Bsingh1}
Basdevant, A-L.,  Singh, A. (2008a) On the speed of a cookie random walk. {\em Prob. Th. Rel. Fields} {\bf 141}, 625--645.

\bibitem{BW}
Benjamini, I., and Wilson, B. (2003) Excited Random Walk. {\em Electron. Comm. Probab.}   86--92.

\bibitem{BHA}  Bhattacharya, R. N. and Waymire, E. (1990)
{\it Stochastic Processes with Applications.} Wiley, New York.

\bibitem{bd90}
Davis, B. (1990)  Reinforced random walk.  {\em Prob. Th. Rel. Fields} {\bf 84}, 203 -- 229.

\bibitem{bd96}
Davis, B. (1996). Weak limits of perturbed random walks and the equation
 $Y\sb t=B\sb t+\alpha\sup\{Y\sb s\colon s\leq t\} +\beta\inf\{Y\sb s\colon s\leq t\}$.
{\em Ann. Probab.} {\bf 24}  2007--2023.

\bibitem{bd99}
Davis, B. (1999)  Brownian motion and  random walk perturbed at extrema.  {\em Prob. Th. Rel. Fields} {\bf 113}, 501 -- 518.

\bibitem{DUR}  Durrett, R. (1996)
{\it Probability: Theory and Examples.} 2nd edition, Wiley, Duxbury press.

\bibitem{WF}
Feller, W. (1971). {\it An Introduction to Probability Theory and Its Applications}, Vol. 2 (second edition). John Wiley.

\bibitem{KosZ}
Kosygina, E., Zerner, M. (2008) Positively and negatively excited random walks on integers, with branching processes. {\it Electr. J. Probab.} {\bf 13}, paper no.~64,  1952--1979.

\bibitem{STH}
Surhone, L.M.,  Tennoe, M.T., and Henssonow, S.F. (Ed.) (2010). {\it Snake Video Game.} Betascript Publishing. 

\bibitem{SV}
Volkov, S. (2003) Excited Random Walks on Trees. {\em Electron. J. Probab.}  {\bf 8}, paper no.~23.

\bibitem{Zern}
Zerner, M. (2005) Multi-excited random walks on integers. {\it Probab. Theory Related Fields} {\bf 133}, 98--122.
\end{thebibliography}
\end{document}